\documentclass[Svgc]{Svgc}

\usepackage{graphicx} % avoid epsfig or earlier such packages
\usepackage{amsfonts}  % many want amsmath extensions

\journalname{Graphs and Combinatorics}
\newcommand{\ds}{\displaystyle}

\begin{document}

\title{Clustering of spectra and fractals of regular graphs}
\author{V. Ejov\inst{1} \and J.A. Filar\inst{1} \and S.K.
Lucas\inst{2} \and P. Zograf\inst{3}}

 \institute{School of
Mathematics and Statistics, University of South Australia,\\
Mawson Lakes SA 5095 \textsc{Australia} \and
 Department of Mathematics and Statistics, James Madison University,\\
 Harrisonburg, VA 22807, \textsc{USA}\and
 St. Petersburg
Department of Steklov Institute of Mathematics, \\
Fontanka 27, St.Petersburg 191023 \textsc{Russia}}

\maketitle

\begin{abstract}
We exhibit a characteristic structure of the class of all regular
graphs of degree $d$ that stems from the spectra of their
adjacency matrices. The structure has a fractal threadlike
appearance. Points with coordinates given by the mean and variance
of the exponentials of graph eigenvalues cluster around a line
segment that we call a {\em filar}. Zooming-in reveals that this
cluster splits into smaller segments (filars) labeled by  the
number of triangles in graphs. Further zooming-in shows that the
smaller filars split into subfilars labelled by the number of
quadrangles in graphs, etc. We call this fractal structure,
discovered in a numerical experiment, a {\em multifilar
structure}. We also provide a mathematical explanation of this
phenomenon based on the Ihara-Selberg trace formula, and compute
the coordinates and slopes of all filars in terms of Bessel
functions of the first kind.
\end{abstract}


\receive{September, 2006}
%\finalreceive{December 22, 2002}

\section{A numerical experiment}

For the sake of simplicity we will pay our attention mainly to
cubic graphs (or, in other words, to regular graphs of degree
$d=3$). This assumption is not restrictive since all our
considerations remain valid for regular graphs of degree $d>3$
with obvious minor modifications. Moreover, in a certain sense
cubic graphs are the generic regular graphs (see e.g. Greenlaw and
Petreschi \cite{Greenlaw})\footnote{Actually, cubic graphs are
generic in a wider sense: {\em any} graph can be made cubic by a
small perturbation that blows up vertices into small circles.}.

So let us consider the set of all regular cubic graphs with $n$ vertices.
They can be conveniently enumerated using the \texttt{GENREG} program of
Markus Meringer \cite{Meringer}. Each graph is completely determined by its
adjacency matrix, which is symmetric. Its spectrum (the set of eigenvalues)
is real and lies on the segment $[-3,3]$. For each graph it can be found
numerically. In the interests of statistical analysis, we might want to
take the means and variances of each set of eigenvalues. However, since
the diagonal entries of the adjacency matrices are zero (graphs contain no
loops), the eigenvalues sum to zero. In order to produce results with
some variation, and originally motivated by solving systems of linear
first order differential equations, we take the exponential of the
eigenvalues before finding their mean and variance. As a final
modification, this time motivated by the authors' interest in Markov
processes, we replace the adjacency matrix $A$ by the related doubly
stochastic matrix $\frac13A$. The theory of Markov chains then states
that the probability of being at the $j$th vertex after a walk of
length $i$ in the graph with each edge equally likely to be chosen is
the $j$th element of the vector $\left(\frac13A\right)^i \mathbf{x}$,
where the $k$th element of the vector $\mathbf{x}$ is the probability
of starting at the $k$th vertex.

\begin{figure}
\includegraphics[width=12.6cm]{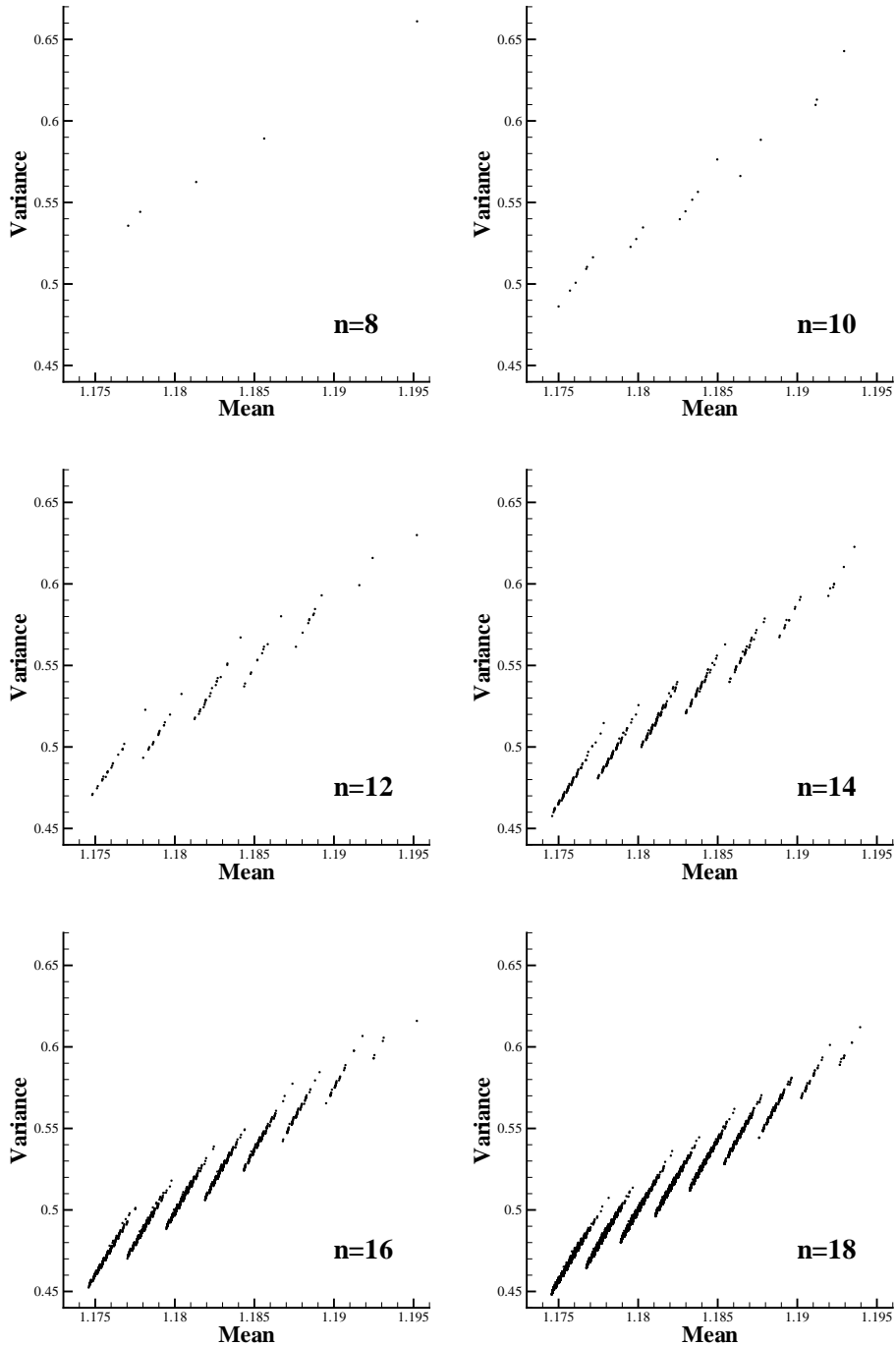}
\caption{Plots of mean versus variance for the exponential of the
eigenvalues of the doubly stochastic matrices associated with all
regular cubic graphs with various numbers of vertices.}
\end{figure}

Summarizing, we apply the following procedure. For a fixed even
$n$ find the adjacency matrices of all regular cubic graphs on $n$
vertices. In each case, divide the adjacency matrix by three, find
its eigenvalues, take their exponential, and then find their mean
and variance. Each cubic graph is then represented by a single dot
on a plot of mean versus variance. Figure 1 shows the results of
applying this procedure with $n=8,10,12,14,16,18$, where the
number of regular cubic graphs in each case is
$5,19,85,509,4060,41301$ respectively. There appears to be a very
definite structure in these plots. In each case the data appear in
distinct clusters that at this scale look like straight line
segments with roughly the same slope and distance separating them.
(In the next section we will derive explicit formulas for these
slopes and distances.) Due to their form, we would like to name
these clusters by ``filars'', whose dictionary meaning is
``threadlike objects''\footnote {This term also recognizes the
second author's (JF's) initial investigation of the above
phenomenon.}.

\begin{figure}
\includegraphics[width=13cm]{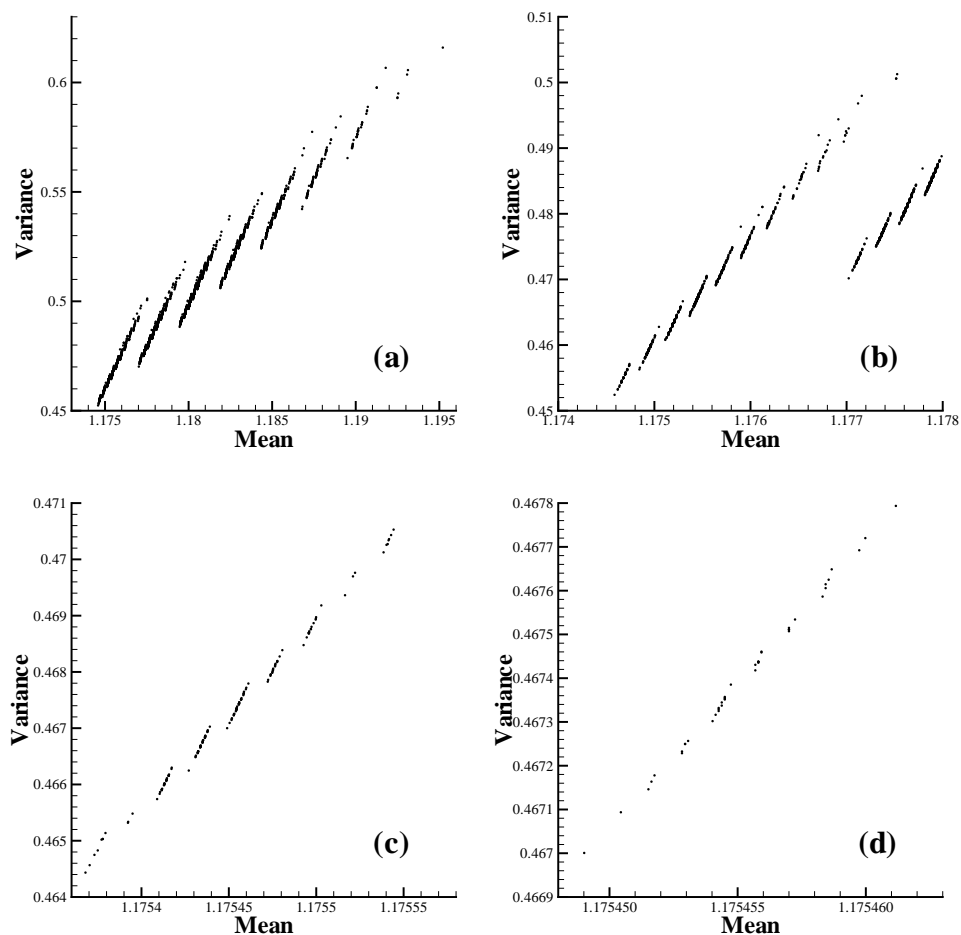}
\caption{Successively zooming in on the $n=14$ plot}
\end{figure}

An even greater level of structure exists within each filar. Figure
2(a) repeats the results for $n=16$, and Figure 2(b) zooms in on the
leftmost filar. We can see that each filar is in fact made up of
smaller clusters of approximately straight line segments, all roughly parallel
and the same distance apart, with a steeper slope than the original one. We
shall call each of these clusters a subfilar, and Figure 2(c) zooms in
on the fourth subfilar from the left in Figure 2(b). The structure
continues in Figure 2(d), which zooms in on the 5th subsubfilar of
Figure 2(c). Since a fractal is defined as a self-similar image, where
the same structure is evident when magnifying one part of the image, we
see that these figures obviously enjoy a fractal structure. The larger the
number of vertices, the more levels of magnification can be undertaken
before the number of data points becomes small enough for the
self-similar structure to be lost. Collectively, we refer to this
phenomenon as the ``multifilar structure'' of cubic graphs (or their
spectra, to be more precise).

Finally, it is worth noting that this behavior is not limited to
cubic graphs. Plots for quartic graphs (every vertex of degree
four) show exactly the same structure. As we will see later, the
Ihara-Selberg trace formula justifies the presence of such a
fractal structure for regular graphs of arbitrary degree $d$.

\setlength{\unitlength}{0.38cm}
\begin{figure}
\begin{center}
\begin{picture}(8,10) % (1)
 \put(3,2){\circle*{0.3}} \put(6,2){\circle*{0.3}}
 \put(2,5){\circle*{0.3}} \put(3,5){\circle*{0.3}}
 \put(4,5){\circle*{0.3}} \put(5,5){\circle*{0.3}}
 \put(6,5){\circle*{0.3}} \put(7,5){\circle*{0.3}}
 \put(3,8){\circle*{0.3}} \put(6,8){\circle*{0.3}}
 \put(3,2){\line(-1,3){1}} \put(3,2){\line(0,1){6}}
 \put(3,2){\line(1,3){1}} \put(6,2){\line(-1,3){1}}
 \put(6,2){\line(0,1){6}} \put(6,2){\line(1,3){1}}
 \put(3,8){\line(-1,-3){1}} \put(3,8){\line(1,-3){1}}
 \put(6,8){\line(-1,-3){1}} \put(6,8){\line(1,-3){1}}
 \put(2,5){\line(1,0){1}} \put(4,5){\line(1,0){1}}
 \put(6,5){\line(1,0){1}} \put(7,8){(1)}
\end{picture}
\begin{picture}(8,10) % (2)
 \put(4,1){\circle*{0.3}} \put(6,2){\circle*{0.3}}
 \put(7,4){\circle*{0.3}} \put(7,6){\circle*{0.3}}
 \put(6,8){\circle*{0.3}} \put(4,9){\circle*{0.3}}
 \put(2,8){\circle*{0.3}} \put(1,6){\circle*{0.3}}
 \put(1,4){\circle*{0.3}} \put(2,2){\circle*{0.3}}
 \put(4,1){\line(2,1){2}} \put(6,2){\line(1,2){1}}
 \put(7,4){\line(0,1){2}} \put(7,6){\line(-1,2){1}}
 \put(6,8){\line(-2,1){2}} \put(4,9){\line(-2,-1){2}}
 \put(2,8){\line(-1,-2){1}} \put(1,6){\line(0,-1){2}}
 \put(1,4){\line(1,-2){1}} \put(2,2){\line(2,-1){2}}
 \put(6,2){\line(1,4){1}} \put(7,4){\line(-1,4){1}}
 \put(2,2){\line(-1,4){1}} \put(1,4){\line(1,4){1}}
 \put(4,1){\line(0,1){8}} \put(7,8){(2)}
\end{picture}
\begin{picture}(8,10) % (3)
 \put(4,1){\circle*{0.3}} \put(6,2){\circle*{0.3}}
 \put(7,4){\circle*{0.3}} \put(7,6){\circle*{0.3}}
 \put(6,8){\circle*{0.3}} \put(4,9){\circle*{0.3}}
 \put(2,8){\circle*{0.3}} \put(1,6){\circle*{0.3}}
 \put(1,4){\circle*{0.3}} \put(2,2){\circle*{0.3}}
 \put(4,1){\line(2,1){2}} \put(6,2){\line(1,2){1}}
 \put(7,4){\line(0,1){2}} \put(7,6){\line(-1,2){1}}
 \put(6,8){\line(-2,1){2}} \put(4,9){\line(-2,-1){2}}
 \put(2,8){\line(-1,-2){1}} \put(1,6){\line(0,-1){2}}
 \put(1,4){\line(1,-2){1}} \put(2,2){\line(2,-1){2}}
 \qbezier(4,1)(2,8)(2,8) \put(2,2){\line(-1,4){1}}
 \put(1,4){\line(3,5){3}} \put(6,2){\line(1,4){1}}
 \put(7,4){\line(-1,4){1}} \put(7,8){(3)}

\end{picture}
\begin{picture}(8,10) % (4)
 \put(4,1){\circle*{0.3}} \put(6,2){\circle*{0.3}}
 \put(7,4){\circle*{0.3}} \put(7,6){\circle*{0.3}}
 \put(6,8){\circle*{0.3}} \put(4,9){\circle*{0.3}}
 \put(2,8){\circle*{0.3}} \put(1,6){\circle*{0.3}}
 \put(1,4){\circle*{0.3}} \put(2,2){\circle*{0.3}}
 \put(4,1){\line(2,1){2}} \put(6,2){\line(1,2){1}}
 \put(7,4){\line(0,1){2}} \put(7,6){\line(-1,2){1}}
 \put(6,8){\line(-2,1){2}} \put(4,9){\line(-2,-1){2}}
 \put(2,8){\line(-1,-2){1}} \put(1,6){\line(0,-1){2}}
 \put(1,4){\line(1,-2){1}} \put(2,2){\line(2,-1){2}}
 \put(4,1){\line(-1,1){3}} \put(1,6){\line(1,1){3}}
 \put(2,2){\line(0,1){6}} \put(6,2){\line(1,4){1}}
 \put(7,4){\line(-1,4){1}} \put(7,8){(4)}
\end{picture}
\begin{picture}(8,10) % (5)
 \put(4,1){\circle*{0.3}} \put(6,2){\circle*{0.3}}
 \put(7,4){\circle*{0.3}} \put(7,6){\circle*{0.3}}
 \put(6,8){\circle*{0.3}} \put(4,9){\circle*{0.3}}
 \put(2,8){\circle*{0.3}} \put(1,6){\circle*{0.3}}
 \put(1,4){\circle*{0.3}} \put(2,2){\circle*{0.3}}
 \put(4,1){\line(2,1){2}} \put(6,2){\line(1,2){1}}
 \put(7,4){\line(0,1){2}} \put(7,6){\line(-1,2){1}}
 \put(6,8){\line(-2,1){2}} \put(4,9){\line(-2,-1){2}}
 \put(2,8){\line(-1,-2){1}} \put(1,6){\line(0,-1){2}}
 \put(1,4){\line(1,-2){1}} \put(2,2){\line(2,-1){2}}
 \put(4,1){\line(-3,5){3}} \put(1,4){\line(3,5){3}}
 \put(2,2){\line(0,1){6}} \put(6,2){\line(1,4){1}}
 \put(7,4){\line(-1,4){1}} \put(7,8){(5)}
\end{picture}
\begin{picture}(8,10) % (6)
 \put(4,1){\circle*{0.3}} \put(6,2){\circle*{0.3}}
 \put(7,4){\circle*{0.3}} \put(7,6){\circle*{0.3}}
 \put(6,8){\circle*{0.3}} \put(4,9){\circle*{0.3}}
 \put(2,8){\circle*{0.3}} \put(1,6){\circle*{0.3}}
 \put(1,4){\circle*{0.3}} \put(2,2){\circle*{0.3}}
 \put(4,1){\line(2,1){2}} \put(6,2){\line(1,2){1}}
 \put(7,4){\line(0,1){2}} \put(7,6){\line(-1,2){1}}
 \put(6,8){\line(-2,1){2}} \put(4,9){\line(-2,-1){2}}
 \put(2,8){\line(-1,-2){1}} \put(1,6){\line(0,-1){2}}
 \put(1,4){\line(1,-2){1}} \put(2,2){\line(2,-1){2}}
 \put(4,1){\line(1,1){3}} \put(1,6){\line(1,1){3}}
 \put(2,2){\line(0,1){6}} \put(6,2){\line(0,1){6}}
 \put(1,4){\line(3,1){6}} \put(7,8){(6)}
\end{picture}
\begin{picture}(8,10) % (7)
 \put(4,1){\circle*{0.3}} \put(6,2){\circle*{0.3}}
 \put(7,4){\circle*{0.3}} \put(7,6){\circle*{0.3}}
 \put(6,8){\circle*{0.3}} \put(4,9){\circle*{0.3}}
 \put(2,8){\circle*{0.3}} \put(1,6){\circle*{0.3}}
 \put(1,4){\circle*{0.3}} \put(2,2){\circle*{0.3}}
 \put(4,1){\line(2,1){2}} \put(6,2){\line(1,2){1}}
 \put(7,4){\line(0,1){2}} \put(7,6){\line(-1,2){1}}
 \put(6,8){\line(-2,1){2}} \put(4,9){\line(-2,-1){2}}
 \put(2,8){\line(-1,-2){1}} \put(1,6){\line(0,-1){2}}
 \put(1,4){\line(1,-2){1}} \put(2,2){\line(2,-1){2}}
 \put(2,8){\line(1,0){4}} \put(1,6){\line(1,0){6}}
 \put(1,4){\line(1,0){6}} \put(2,2){\line(1,0){4}}
 \put(4,1){\line(0,1){8}} \put(7,8){(7)}
\end{picture}
\begin{picture}(8,10) % (8)
 \put(4,1){\circle*{0.3}} \put(6,2){\circle*{0.3}}
 \put(7,4){\circle*{0.3}} \put(7,6){\circle*{0.3}}
 \put(6,8){\circle*{0.3}} \put(4,9){\circle*{0.3}}
 \put(2,8){\circle*{0.3}} \put(1,6){\circle*{0.3}}
 \put(1,4){\circle*{0.3}} \put(2,2){\circle*{0.3}}
 \put(4,1){\line(2,1){2}} \put(6,2){\line(1,2){1}}
 \put(7,4){\line(0,1){2}} \put(7,6){\line(-1,2){1}}
 \put(6,8){\line(-2,1){2}} \put(4,9){\line(-2,-1){2}}
 \put(2,8){\line(-1,-2){1}} \put(1,6){\line(0,-1){2}}
 \put(1,4){\line(1,-2){1}} \put(2,2){\line(2,-1){2}}
 \put(2,8){\line(1,0){4}} \put(1,6){\line(1,0){6}}
 \put(2,2){\line(5,2){5}} \put(1,4){\line(5,-2){5}}
 \put(4,1){\line(0,1){8}} \put(7,8){(8)}
\end{picture}
\begin{picture}(8,10) % (9)
 \put(4,1){\circle*{0.3}} \put(6,2){\circle*{0.3}}
 \put(7,4){\circle*{0.3}} \put(7,6){\circle*{0.3}}
 \put(6,8){\circle*{0.3}} \put(4,9){\circle*{0.3}}
 \put(2,8){\circle*{0.3}} \put(1,6){\circle*{0.3}}
 \put(1,4){\circle*{0.3}} \put(2,2){\circle*{0.3}}
 \put(4,1){\line(2,1){2}} \put(6,2){\line(1,2){1}}
 \put(7,4){\line(0,1){2}} \put(7,6){\line(-1,2){1}}
 \put(6,8){\line(-2,1){2}} \put(4,9){\line(-2,-1){2}}
 \put(2,8){\line(-1,-2){1}} \put(1,6){\line(0,-1){2}}
 \put(1,4){\line(1,-2){1}} \put(2,2){\line(2,-1){2}}
 \put(2,8){\line(1,0){4}} \put(2,2){\line(1,0){4}}
 \put(1,4){\line(3,1){6}} \put(1,6){\line(3,-1){6}}
 \put(4,1){\line(0,1){8}} \put(7,8){(9)}
\end{picture}
\begin{picture}(8,10) % (10)
 \put(4,1){\circle*{0.3}} \put(6,2){\circle*{0.3}}
 \put(7,4){\circle*{0.3}} \put(7,6){\circle*{0.3}}
 \put(6,8){\circle*{0.3}} \put(4,9){\circle*{0.3}}
 \put(2,8){\circle*{0.3}} \put(1,6){\circle*{0.3}}
 \put(1,4){\circle*{0.3}} \put(2,2){\circle*{0.3}}
 \put(4,1){\line(2,1){2}} \put(6,2){\line(1,2){1}}
 \put(7,4){\line(0,1){2}} \put(7,6){\line(-1,2){1}}
 \put(6,8){\line(-2,1){2}} \put(4,9){\line(-2,-1){2}}
 \put(2,8){\line(-1,-2){1}} \put(1,6){\line(0,-1){2}}
 \put(1,4){\line(1,-2){1}} \put(2,2){\line(2,-1){2}}
 \put(4,1){\line(0,1){8}} \put(1,4){\line(1,4){1}}
 \put(7,4){\line(-1,4){1}} \put(2,2){\line(5,4){5}}
 \put(6,2){\line(-5,4){5}} \put(7,8){(10)}
\end{picture}
\begin{picture}(8,10) % (11)
 \put(4,1){\circle*{0.3}} \put(6,2){\circle*{0.3}}
 \put(7,4){\circle*{0.3}} \put(7,6){\circle*{0.3}}
 \put(6,8){\circle*{0.3}} \put(4,9){\circle*{0.3}}
 \put(2,8){\circle*{0.3}} \put(1,6){\circle*{0.3}}
 \put(1,4){\circle*{0.3}} \put(2,2){\circle*{0.3}}
 \put(4,1){\line(2,1){2}} \put(6,2){\line(1,2){1}}
 \put(7,4){\line(0,1){2}} \put(7,6){\line(-1,2){1}}
 \put(6,8){\line(-2,1){2}} \put(4,9){\line(-2,-1){2}}
 \put(2,8){\line(-1,-2){1}} \put(1,6){\line(0,-1){2}}
 \put(1,4){\line(1,-2){1}} \put(2,2){\line(2,-1){2}}
 \put(4,1){\line(0,1){8}} \put(2,2){\line(0,1){6}}
 \put(1,4){\line(1,0){6}} \put(1,6){\line(5,2){5}}
 \put(6,2){\line(1,4){1}} \put(7,8){(11)}
\end{picture}
\begin{picture}(8,10) % (12)
 \put(4,1){\circle*{0.3}} \put(6,2){\circle*{0.3}}
 \put(7,4){\circle*{0.3}} \put(7,6){\circle*{0.3}}
 \put(6,8){\circle*{0.3}} \put(4,9){\circle*{0.3}}
 \put(2,8){\circle*{0.3}} \put(1,6){\circle*{0.3}}
 \put(1,4){\circle*{0.3}} \put(2,2){\circle*{0.3}}
 \put(4,1){\line(2,1){2}} \put(6,2){\line(1,2){1}}
 \put(7,4){\line(0,1){2}} \put(7,6){\line(-1,2){1}}
 \put(6,8){\line(-2,1){2}} \put(4,9){\line(-2,-1){2}}
 \put(2,8){\line(-1,-2){1}} \put(1,6){\line(0,-1){2}}
 \put(1,4){\line(1,-2){1}} \put(2,2){\line(2,-1){2}}
 \put(4,1){\line(0,1){8}} \put(2,2){\line(-1,4){1}}
 \put(6,2){\line(1,4){1}} \put(2,8){\line(1,0){4}}
 \put(1,4){\line(1,0){6}} \put(7,8){(12)}
\end{picture}
\begin{picture}(8,10) % (13)
 \put(4,1){\circle*{0.3}} \put(6,2){\circle*{0.3}}
 \put(7,4){\circle*{0.3}} \put(7,6){\circle*{0.3}}
 \put(6,8){\circle*{0.3}} \put(4,9){\circle*{0.3}}
 \put(2,8){\circle*{0.3}} \put(1,6){\circle*{0.3}}
 \put(1,4){\circle*{0.3}} \put(2,2){\circle*{0.3}}
 \put(4,1){\line(2,1){2}} \put(6,2){\line(1,2){1}}
 \put(7,4){\line(0,1){2}} \put(7,6){\line(-1,2){1}}
 \put(6,8){\line(-2,1){2}} \put(4,9){\line(-2,-1){2}}
 \put(2,8){\line(-1,-2){1}} \put(1,6){\line(0,-1){2}}
 \put(1,4){\line(1,-2){1}} \put(2,2){\line(2,-1){2}}
 \put(4,1){\line(0,1){8}} \put(1,4){\line(1,0){6}}
 \put(2,2){\line(5,4){5}} \put(1,6){\line(5,-4){5}}
 \put(2,8){\line(1,0){4}} \put(7,8){(13)}
\end{picture}
\begin{picture}(8,10) % (14)
 \put(4,1){\circle*{0.3}} \put(6,2){\circle*{0.3}}
 \put(7,4){\circle*{0.3}} \put(7,6){\circle*{0.3}}
 \put(6,8){\circle*{0.3}} \put(4,9){\circle*{0.3}}
 \put(2,8){\circle*{0.3}} \put(1,6){\circle*{0.3}}
 \put(1,4){\circle*{0.3}} \put(2,2){\circle*{0.3}}
 \put(4,1){\line(2,1){2}} \put(6,2){\line(1,2){1}}
 \put(7,4){\line(0,1){2}} \put(7,6){\line(-1,2){1}}
 \put(6,8){\line(-2,1){2}} \put(4,9){\line(-2,-1){2}}
 \put(2,8){\line(-1,-2){1}} \put(1,6){\line(0,-1){2}}
 \put(1,4){\line(1,-2){1}} \put(2,2){\line(2,-1){2}}
 \put(4,1){\line(0,1){8}} \put(2,2){\line(5,2){5}}
 \put(1,4){\line(5,-2){5}} \put(1,6){\line(5,2){5}}
 \put(2,8){\line(5,-2){5}} \put(7,8){(14)}
\end{picture}
\begin{picture}(8,10) % (15)
 \put(4,1){\circle*{0.3}} \put(6,2){\circle*{0.3}}
 \put(7,4){\circle*{0.3}} \put(7,6){\circle*{0.3}}
 \put(6,8){\circle*{0.3}} \put(4,9){\circle*{0.3}}
 \put(2,8){\circle*{0.3}} \put(1,6){\circle*{0.3}}
 \put(1,4){\circle*{0.3}} \put(2,2){\circle*{0.3}}
 \put(4,1){\line(2,1){2}} \put(6,2){\line(1,2){1}}
 \put(7,4){\line(0,1){2}} \put(7,6){\line(-1,2){1}}
 \put(6,8){\line(-2,1){2}} \put(4,9){\line(-2,-1){2}}
 \put(2,8){\line(-1,-2){1}} \put(1,6){\line(0,-1){2}}
 \put(1,4){\line(1,-2){1}} \put(2,2){\line(2,-1){2}}
 \put(2,2){\line(0,1){6}} \put(4,1){\line(0,1){8}}
 \put(6,2){\line(0,1){6}} \put(1,4){\line(1,0){6}}
 \put(1,6){\line(1,0){6}} \put(7,8){(15)}
\end{picture}
\begin{picture}(8,10) % (16)
 \put(4,1){\circle*{0.3}} \put(6,2){\circle*{0.3}}
 \put(7,4){\circle*{0.3}} \put(7,6){\circle*{0.3}}
 \put(6,8){\circle*{0.3}} \put(4,9){\circle*{0.3}}
 \put(2,8){\circle*{0.3}} \put(1,6){\circle*{0.3}}
 \put(1,4){\circle*{0.3}} \put(2,2){\circle*{0.3}}
 \put(4,1){\line(2,1){2}} \put(6,2){\line(1,2){1}}
 \put(7,4){\line(0,1){2}} \put(7,6){\line(-1,2){1}}
 \put(6,8){\line(-2,1){2}} \put(4,9){\line(-2,-1){2}}
 \put(2,8){\line(-1,-2){1}} \put(1,6){\line(0,-1){2}}
 \put(1,4){\line(1,-2){1}} \put(2,2){\line(2,-1){2}}
 \put(2,2){\line(2,3){4}} \put(4,1){\line(0,1){8}}
 \put(6,2){\line(-2,3){4}} \put(1,4){\line(1,0){6}}
 \put(1,6){\line(1,0){6}} \put(7,8){(16)}
\end{picture}
\begin{picture}(8,10) % (17)
 \put(4,1){\circle*{0.3}} \put(6,2){\circle*{0.3}}
 \put(7,4){\circle*{0.3}} \put(7,6){\circle*{0.3}}
 \put(6,8){\circle*{0.3}} \put(4,9){\circle*{0.3}}
 \put(2,8){\circle*{0.3}} \put(1,6){\circle*{0.3}}
 \put(1,4){\circle*{0.3}} \put(2,2){\circle*{0.3}}
 \put(4,1){\line(2,1){2}} \put(6,2){\line(1,2){1}}
 \put(7,4){\line(0,1){2}} \put(7,6){\line(-1,2){1}}
 \put(6,8){\line(-2,1){2}} \put(4,9){\line(-2,-1){2}}
 \put(2,8){\line(-1,-2){1}} \put(1,6){\line(0,-1){2}}
 \put(1,4){\line(1,-2){1}} \put(2,2){\line(2,-1){2}}
 \put(2,2){\line(2,3){4}} \put(4,1){\line(0,1){8}}
 \put(6,2){\line(-2,3){4}} \put(1,4){\line(3,1){6}}
 \put(1,6){\line(3,-1){6}} \put(7,8){(17)}
\end{picture}
\begin{picture}(8,10) % (18)
 \put(4,1){\circle*{0.3}} \put(6,2){\circle*{0.3}}
 \put(7,4){\circle*{0.3}} \put(7,6){\circle*{0.3}}
 \put(6,8){\circle*{0.3}} \put(4,9){\circle*{0.3}}
 \put(2,8){\circle*{0.3}} \put(1,6){\circle*{0.3}}
 \put(1,4){\circle*{0.3}} \put(2,2){\circle*{0.3}}
 \put(4,1){\line(2,1){2}} \put(6,2){\line(1,2){1}}
 \put(7,4){\line(0,1){2}} \put(7,6){\line(-1,2){1}}
 \put(6,8){\line(-2,1){2}} \put(4,9){\line(-2,-1){2}}
 \put(2,8){\line(-1,-2){1}} \put(1,6){\line(0,-1){2}}
 \put(1,4){\line(1,-2){1}} \put(2,2){\line(2,-1){2}}
 \put(2,2){\line(5,4){5}} \put(1,4){\line(5,4){5}}
 \put(1,6){\line(5,-4){5}} \put(2,8){\line(5,-4){5}}
 \put(4,1){\line(0,1){8}} \put(7,8){(18)}
\end{picture}
\begin{picture}(8,10) % (19)
 \put(1,6){\circle*{0.3}} \put(2,2){\circle*{0.3}}
 \put(2,4){\circle*{0.3}} \put(2,6){\circle*{0.3}}
 \put(4,8){\circle*{0.3}} \put(4,9){\circle*{0.3}}
 \put(6,2){\circle*{0.3}} \put(6,4){\circle*{0.3}}
 \put(6,6){\circle*{0.3}} \put(7,6){\circle*{0.3}}
 \put(2,2){\line(0,1){2}} \put(6,2){\line(0,1){2}}
 \put(1,6){\line(1,0){1}} \put(6,6){\line(1,0){1}}
 \put(4,8){\line(0,1){1}} \put(6,2){\line(1,4){1}}
 \put(7,6){\line(-1,1){3}} \put(4,9){\line(-1,-1){3}}
 \put(1,6){\line(1,-4){1}} \put(2,2){\line(1,0){4}}
 \put(2,6){\line(1,0){4}} \put(2,6){\line(2,-1){4}}
 \put(2,4){\line(1,2){2}} \put(2,4){\line(2,1){4}}
 \put(6,4){\line(-1,2){2}} \put(7,8){(19)}
\end{picture}
\end{center}
\caption{All cubic graphs with ten vertices, as output by
\texttt{GENREG}} \label{fig:cubics}
\end{figure}
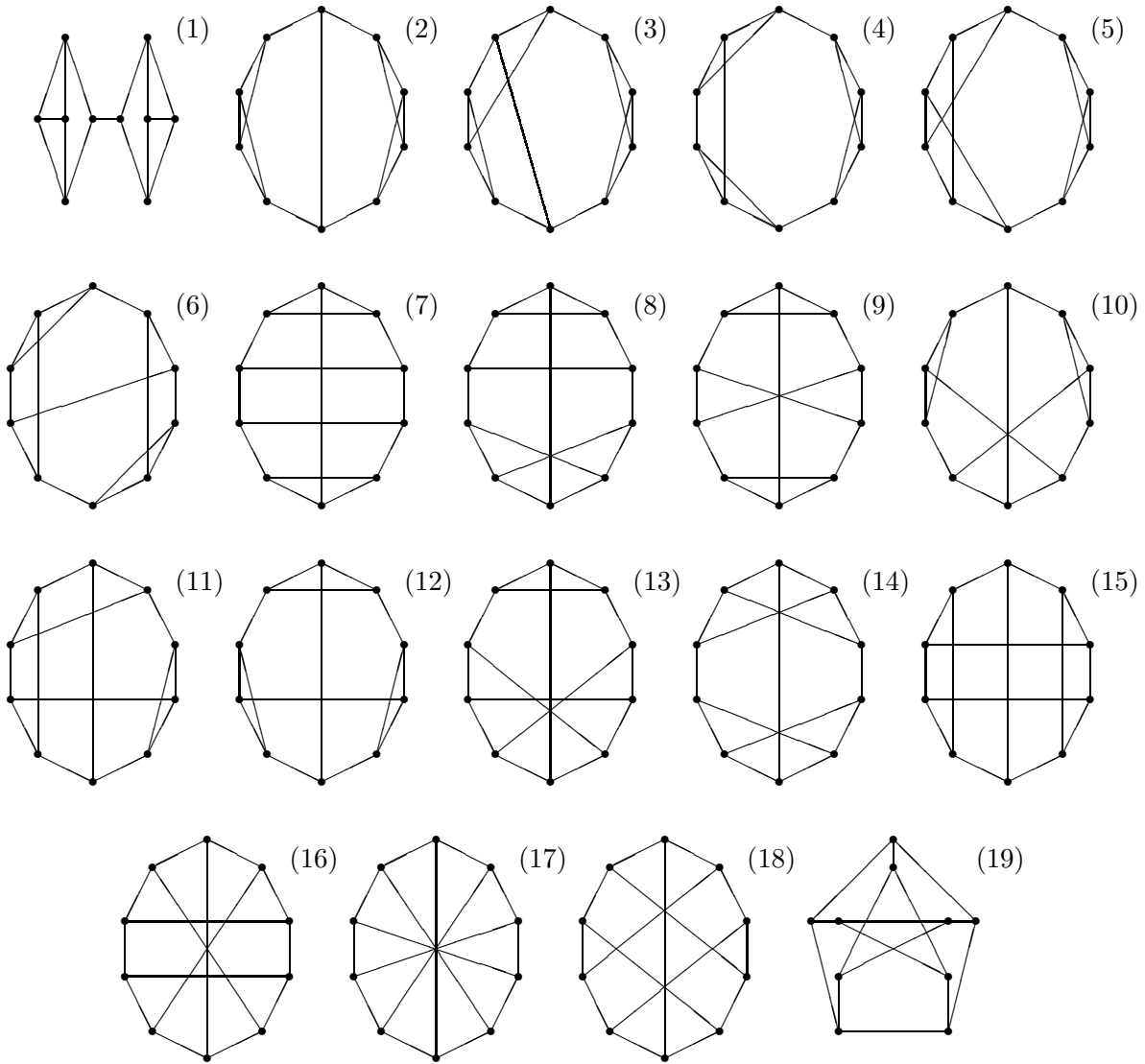

Having -- to the best of our knowledge -- for the first time discovered
this property of regular graphs, the aim of the following section is to
theoretically explain the (multi)filar structure.

\section{Theoretical justification}

\subsection{Ten vertex cubic graphs in detail}

\begin{figure}
\centerline{\includegraphics[width=13cm]{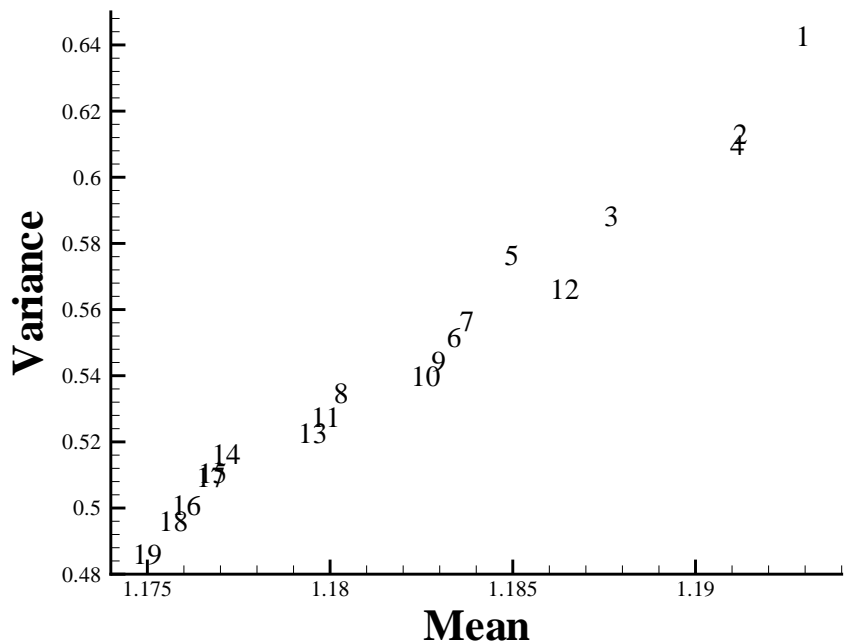}} \caption{A
reproduction of Figure 1($n=10$) with labels replacing data points
related to associated with graphs in Figure 3}
\end{figure}

Before we resort to the theory based on the Ihara-Selberg trace
formula \cite{Ihara,Ahumada}, it is instructive to consider the
case $n=10$ in detail. Figure 3 shows all 19 regular cubic graphs
with ten vertices, labelled in the order produced by
\texttt{GENREG}, and Figure 4 repeats Figure 1(b), but with labels
on the data points indicating graph number. We can see a pattern
in terms of which graphs are in each filar. Graphs
19,18,16,17,15,14 in the first (leftmost) filar have no subcycles
of length 3, which we shall call triangles from now on. Graphs
13,11,8 in the second filar have exactly one triangle, graphs
10,9,6,7,5 in the third filar have exactly two triangles, graphs
12,3 in the fourth filar have exactly three triangles, and graphs
4,2,1 in the fifth filar have exactly four triangles. Focusing on
the first filar, graphs 19,18,16,17,15,14 have exactly 0,2,3,5,5,6
subcycles of length 4 respectively. Special emphasis should be
given to the data points for graphs 17 and 15. They are extremely
close together, and have the same number of subcycles of lengths 3
and 4 (0 and 5 respectively). They only vary in number of
subcycles of length 5, numbering 0 and 2. Similarly, graphs 4 and
2 both have four 3-subcycles, two 4-subcycles, and only vary in
the number of 5-subcycles (two and four). These observations
suggest that membership in a filar structure is related to the
number of subcycles of various lengths in the original graph --
although this does not explain why filars approximate straight
lines.

\subsection{Explicit formulas for the mean and variance}

To obtain qualitative and quantitave justification of the
phenomenon in question we bring in a very explicit version of the
Ihara-Selberg trace formula that is due to P. Mn\"ev, cf.
\cite{Mnev}, Formula (30). A general form of the Ihara-Selberg
trace formula, as well as precise definitions can be found in
Appendix. For any regular graph $G$ of degree $d=q+1$ on $n$
vertices we have

\begin{equation}
\frac{1}{n}\sum_{i=1}^n e^{t\lambda_i}
=\frac{q+1}{2\pi}\int_{-2\sqrt{q}}^{2\sqrt{q}}
e^{st}\frac{\sqrt {4q-s^2}}{(q+1)^2-s^2}ds
+\frac{1}{n}\sum_\gamma
\sum_{k=1}^{\infty}\frac{\ell(\gamma)}{2^{k\ell(\gamma)/2}}
I_{k\ell(\gamma)}(2\sqrt{q}t).\label{Mnev}
\end{equation}

Here $\{\lambda_1,\dots,\lambda_n\}$ is the spectrum of the
adjacency matrix of $G$, $\gamma$ runs over the set of all
(oriented) primitive closed geodesics\footnote {In the context of
graphs, a {\em closed geodesic} is an oriented closed path of
minimal length in its free homotopy class. A closed geodesic is
called {\em primitive} if it is not a multiple of a shorter
geodesic. Closed geodesics of length 3, 4 and 5 are 3-, 4- and
5-cycles in graphs respectively, whereas closed geodesics of
length greater than 5 may have self-intersections. See Appendix
for details.} in $G$, $\ell(\gamma)$ is the length of $\gamma$,
and $I_m(z)$ is the standard notation for the Bessel function of
the first kind:

$$I_m(z)=\sum_{r=0}^\infty\frac{(z/2)^{n+2r}}{r!(n+r)!}.$$

All lengths $\ell(\gamma)$ are integers greater than or equal to
3. Let us denote by $m_\ell$ the number of {\em non-oriented}
primitive closed geodesics of length $\ell$ in the graph $G$. The
numbers $m_\ell$ are called the {\em multiplicities} of the {\em
length spectrum} of the graph $G$, that is, the set of lengths of
non-oriented primitive closed geodesics in $G$ (the set
$\{m_3,m_4,\dots\}$ describes the length spectrum of the graph in
a unique and convenient way). The multiplicities $m_i$ are
uniquely determined by the eigenvalues of $G$ and are given by
explicit formulas (cf. Formula (34) in
\cite{Mnev})\footnote{Clearly, the length spectrum also determines
the eigenvalue spectrum uniquely.}. For instance, we have
$$n_3=\frac{1}{6}\sum_{i=1}^n\lambda_i^3,\quad
n_4=\frac{1}{8}\left(\sum_{i=1}^n\lambda_i^4-n(q+1)(2q+1)\right),$$
etc.

Now we rewrite Formula (\ref{Mnev}) in terms of the multiplicities $n_\ell$
(since we consider in detail only the case of cubic graphs, we also put $q=2$):

\begin{equation}
\frac{1}{n}\sum_{i=1}^n e^{t\lambda_i}
=J(t)+\frac{2}{n}\sum_{\ell=3}^{\infty}\ell m_\ell F_\ell(t),
\label{cubic}\end{equation}

where

$$J(t) = \frac{3}{2\pi}\int\limits^{2\sqrt{2}}_{-2\sqrt{2}}
e^{st}\frac{\sqrt {8-s^2}}{9-s^2}ds,$$

and

$$F_\ell(t)
=\sum_{k=1}^\infty \frac{I_{k\ell}(2\sqrt{2}t)}{2^{k\ell/2}}.$$
Note that the factor of 2 at the sum appears because we forget
about the orientation of geodesics and have to count each one of
them twice. The latter series converges very fast because of the
following well-known asymptotic behavior of $I_m(z)$:

\begin{equation}
I_m(z) \approx \frac{1}{m!} \left(\frac{z}{2}\right)^m
\label{asympt}\end{equation}

as $m\rightarrow\infty$ and $0<z \ll\sqrt{m + 1}$; see e.g. \cite{WW}.

The closed form expressions for the mean $\mu$ and the variance
$\sigma$ can now be easily extracted from (\ref{cubic}).
Precisely, we have
\begin{equation}\label{mean}
\mu = \frac{1}{n}\sum_{i=1}^n e^{\lambda_i/3}
=J(1/3)+\frac{2}{n}\sum_{\ell=3}^{\infty}\ell m_\ell F_\ell(1/3),
\end{equation}

and\footnote{It should be mentioned that the plots on Figures 1
and 2 are build using the {\em unbiased} variance
$s^2_{n-1}=\frac{1}{n-1} \sum_{i=1}^n \left(
e^{\frac{\lambda_i}{3}}-\mu\right)^2$. In order to make formulas
simpler, we consider here the variance $s^2_{n}=\frac{1}{n}
\sum_{i=1}^n \left( e^{\frac{\lambda_i}{3}}-\mu\right)^2$. The
difference is insignificant, especially for large enough $n$.}

\begin{equation}
\begin{array}{ll} \sigma \! & \ds
= \frac{1}{n} \sum_{i=1}^n \left( e^{\frac{\lambda_i}{3}}-\mu\right)^2
= \frac{1}{n}\sum_{i=1}^n \left(e^{\frac{2\lambda_i}{3}}\right) -\mu^2\\ & \ds
= J(2/3)+\frac{2}{n} \sum\limits^n_{\ell=3} \ell n_\ell F_\ell(2/3) -\mu^2.
\end{array}
\end{equation}

Substituting (\ref{mean}) into the last formula and neglecting
quadratic terms in $F_\ell$ that are small in view of
(\ref{asympt}), we get
\begin{equation}\label{variance}
\sigma \approx \left(J(2/3)-J(1/3)^2\right)+\frac{2}{n} \sum
\limits^n_{\ell=1} \ell m_\ell
\left(F_\ell(2/3)-2J(1/3)F_\ell(1/3)\right).
\end{equation}

Now we are set for explicitly describing the positions of filars.

\subsection{Coordinates and slopes of filars}

To start with, let us note that the function $F_\ell(t)$ is
positive for $t > 0$ and decreases very rapidly when $\ell$ grows
and $t$ remains fixed. It is easy to check that
$F_\ell(2/3)-2J(1/3)F_\ell(1/3)$ is also positive for any positive
integer $\ell\geq 3$ and decreases very fast in $\ell$. Therefore,
for any $n$ all the points with coordinates $(\mu,\sigma)$,
corresponding to cubic graphs on $n$ vertices, lie above and to
the right of the initial point $(J(1/3),J(2/3)-J(1/3)^2)\approx
(1.17455, 0.4217)$ in the mean vs. variance plane. Moreover, we
see that $\ell=3$ gives the leading terms in both sums in
(\ref{mean}) and (\ref{variance}). This means that the points
$(\mu,\sigma)$ accumulate just above the line parametrically
described by equations
$$\begin{array}{ll} x \! & \ds
= J(1/3)+tF_3(1/3)\approx 1.17455+0.00653t,\\ \\
y \!& \ds
= \left(J(2/3)-J(1/3)^2\right)+t\left(F_3(2/3)-2J(1/3)F_3(1/3)\right)
\approx 0.4217+0.0462t,
\end{array}$$
where $x$ is the mean and $y$ is the variance coordinates respectively. Note
that the slope of this line is approximately 7.079.

Another byproduct of the above considerations is a necessary and
sufficient condition for two graphs $G^{(1)}$ and $G^{(2)}$ to
belong to the same filar: this happens if and only if the
multiplicities $n_3^{(1)}$ and $n_3^{(2)}$ are equal, or,
equivalently, $G^{(1)}$ and $G^{(2)}$ have equal number of
triangles. The lines these filars approximate
are given by parametric equations
\vspace{0.1in}
$$\begin{array}{ll} x \!& \ds
= J(1/3)+6n_3F_3(1/3)/n+tF_4(1/3),\\\\\vspace{0.1in}
y \!& \ds
= \left(J(2/3)-J(1/3)^2\right)
+6n_3\left(F_3(2/3)-2J(1/3)F_3(1/3)\right)/n \\
& \ds +\; t\left(F_4(2/3)-2J(1/3)F_4(1/3)\right).
\end{array}$$
The horizontal distance between two filars that contain the points
corresponding to $G^{(1)}$ and $G^{(2)}$ is proportional to
$n_3^{(1)} - n_3^{(2)}$ and is approximately equal to
$$6\left(F_3(1/3)-F_4(1/3)
\frac{F_3(2/3)-2J(1/3)F_3(1/3)}{F_4(2/3)-2J(1/3)F_4(1/3)}\right)
\frac{n_3^{(1)} - n_3^{(2)}}{n}.$$ For $n=12,14,16,18$ the
approximate horizontal distances between the neighboring filars
are 0.00181, 0.00155, 0.00136, 0.00121 respectively. As one can
easily see on Figure 1, filars actually get closer to each other
as $n$ gets larger in proportion with $1/n$. However, the slope of
filars is independent of $n$ and is equal to
$F_4(2/3)/F_4(1/3)-2J(1/3)\approx 15.89$. All the above agree
perfectly with the numerical data plotted in Figures 1 and 2.

Each filar splits into subfilars labelled by the number $n_4$ of
quadrangles in the corresponding graphs. These subfilars
approximate line segments of slope
$F_5(2/3)/F_5(1/3)-2J(1/3)\approx 33.36$. The horizontal distance
between subfilars is measured by increments of
$$\frac{8}{n}\left(F_4(1/3)-F_5(1/3)
\frac{F_4(2/3)-2J(1/3)F_4(1/3)}{F_5(2/3)-2J(1/3)F_5(1/3)}\right).$$
One can pursue this kind of analysis for subfilars of any level (or
depth).

\begin{figure}
\begin{center} \hspace*{1cm}
\begin{picture}(14.5,7)
 \put(4,2){\circle*{0.3}} \put(5,3){\circle*{0.3}}
 \put(4,4){\circle*{0.3}} \put(3,3){\circle*{0.3}}
 \put(3,3){\line(-1,0){1}} \put(4,2){\line(0,1){2}}
 \put(4,2){\line(1,1){1}} \put(4,2){\line(-1,1){1}}
 \put(4,4){\line(1,-1){1}} \put(4,4){\line(-1,-1){1}}
 \put(7,2){\circle*{0.3}} \put(8,3){\circle*{0.3}}
 \put(7,4){\circle*{0.3}} \put(6,3){\circle*{0.3}}
 \put(6,3){\line(-1,0){1}} \put(7,2){\line(0,1){2}}
 \put(7,2){\line(1,1){1}} \put(7,2){\line(-1,1){1}}
 \put(7,4){\line(1,-1){1}} \put(7,4){\line(-1,-1){1}}
 \put(10,2){\circle*{0.3}} \put(11,3){\circle*{0.3}}
 \put(10,4){\circle*{0.3}} \put(9,3){\circle*{0.3}}
 \put(9,3){\line(-1,0){1}} \put(10,2){\line(0,1){2}}
 \put(10,2){\line(1,1){1}} \put(10,2){\line(-1,1){1}}
 \put(10,4){\line(1,-1){1}} \put(10,4){\line(-1,-1){1}}
 \put(11,3){\line(1,0){1}} \put(12.5,5){(a)}
 \put(12.5,3){...} \put(0.7,3){...}
\end{picture} \quad
\begin{picture}(8.5,7)
 \put(3,2){\circle*{0.3}} \put(4,3){\circle*{0.3}}
 \put(3,4){\circle*{0.3}} \put(2,3){\circle*{0.3}}
 \put(5,3){\circle*{0.3}} \put(2,3){\line(1,1){1}}
 \put(2,3){\line(1,-1){1}} \put(2,3){\line(1,0){2}}
 \put(4,3){\line(-1,-1){1}} \put(4,3){\line(-1,1){1}}
 \put(5,3){\line(-2,1){2}} \put(5,3){\line(-2,-1){2}}
 \put(5,3){\line(1,0){1}} \put(6.5,3){...} \put(6.5,5){(b)}
\end{picture}
\end{center}
\caption{A string of diamonds (a) and a clasp (b)} \label{diamond}
\end{figure}
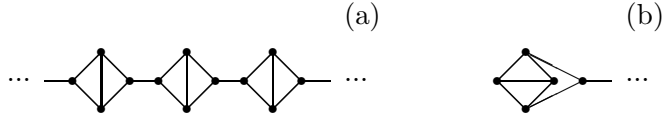

Finally, a comment is in order. The first one concerns the number
of filars, or, equivalently, how many distinct values the
multiplicity $m_3$ (= the number of triangles) can attain for a
regular graph on $n$ vertices. We sketch an argument that gives an
upper bound for $m_3$ in terms of $n$. It can be shown that the
maximal number of triangles is achieved in planar graphs. By the
Euler characteristic formula we have $v-e+f=2$, where $v,e,f$ are
the numbers of vertices, edges and faces of a planar graph. Denote
by $f_k$ the number of its $k$-gonal faces. Then, for a
$d$-regular graph on $n$ vertices $v=n,\;e=nd/2,\; f=\sum f_k$.
Substituting these expressions into the Euler characteristic
formula, we get
$$\left(1-\frac{d}{2}\right)n + \sum_{k\geq3}f_k=2.$$
Say, for $d=3$ a careful analysis of this formula shows that
$m_3=f_3\leq 2[n/4]$ with the exception of $K_4$, the complete
graph on 4 vertices. This upper bound is sharp. For $n$ divisible
by 4 it is achieved by looping a \textit{string of diamonds} shown
in Figure \ref{diamond}(a)(cf. \cite{Greenlaw}, \cite{Korfhage}).
When $n\equiv 2 \;\mbox{mod}\; 4$, we need to attach a
\textit{clasp} on either end of a string of diamonds, as shown in
Figure \ref{diamond}(b). In fact, for cubic graphs on $n\geq 8$
vertices the multiplicity $m_3$ can be any number between 0 and
$2[n/4]$; the corresponding examples can also be easily
constructed (e.g. $m_3=0$ for bipartite graphs).

%The second comment is of a more sophisticated nature. Here we
%tacitly assumed that the contribution from longer geodesics to the
%trace formula (\ref{Mnev}) is always small compared to the
%contribution from shorter ones. Actually, this is a rather
%delicate issue that we hope to address in detail elsewhere.

\bigskip\bigskip
\appendix{{\bf Appendix.} {\it The Ihara-Selberg trace formula}}
\bigskip

The famous Selberg trace formula relates the eigenvalue spectrum of
the Laplace operator on a hyperbolic surface to its length spectrum --
the collection of lengths of closed geodesics counted with multiplicities.
An immediate consequence of the Selberg trace formula is that the
eigenvalue spectrum and the length spectrum uniquely determine one another.
A similar result is valid for regular graphs \cite{Ihara}, \cite{Ahumada}.
To formulate it precisely we need to introduce some terminology.

We consider oriented closed paths in graphs up to cyclic permutations
of vertices. An elementary homotopy is a transformation of a closed path
of the form
$$ (v_1,\ldots,v_j,\ldots,v_n,v_1)
\mapsto (v_1,\ldots,v_j,v',v_j\ldots,v_n,v_1),$$ where $v_j$ and
$v'$ are adjacent vertices. Two closed paths are called (freely)
homotopic if one can be transformed into another by a sequence of
elementary homotopies or their inverses. The unique shortest
representative in a (free) homotopy class of closed paths is
called a {\em closed geodesic}. The {\em length} $\ell(\gamma)$ of
a closed geodesic $\gamma$ is the number of edges it contains;
$\gamma$ is called {\em primitive} if it is not a power of a
shorter geodesic.

Now let $G$ be a regular graph of degree $d=q+1$ on $n$ vertices. Denote by
$A$ its adjacency matrix, and let $\{\lambda_1>\lambda_2\geq\dots\lambda_n\}$
be the spectrum of $A$. Note that $\lambda_1=q+1$ and $|\lambda_i|\leq q+1$.
The following result can be found in \cite{Ahumada}:

\begin{theorem} Let $h:\mathbb{Z} \rightarrow
\mathbb{C}$ be a sequence of complex numbers such that $h(n)=h(-n)$ for
all $n\in\mathbb{Z}$ and
$$\sum_{n=1}^\infty |h(n)| q^{n/2} < \infty.$$
Put $\hat{h}(z)=\sum_{n=-\infty}^\infty h(n)z^{-n}$, the discrete Fourier transform
of $h(n)$. Then
\begin{equation}
\sum_{i=1}^n \hat{h}(z_i)=\frac{nq}{2\pi\sqrt{-1}}
\oint_{|z|=1}\hat{h}(z)\frac{1-z^2}{q-z^2}\,\frac{dz}{z}
+ \sum_{\gamma}\sum_{k=1}^\infty \frac{\ell(\gamma)}{q^{n\ell(\gamma)/2}}h(k\ell(\gamma)),
\label{Ahumada}\end{equation}
where $z_i$ is related to $\lambda_i$ by the equation $\lambda_i=\sqrt{q}(z_i+z_i^{-1})$,
and $\gamma$ runs over the set of all primitive closed geodesics in $G$.
\end{theorem}

Formula (\ref{Mnev}) follows from (\ref{Ahumada}) if we take $h(n)=I_n(2\sqrt{q}t)$.
Then we have $\hat{h}(z)=e^{t\sqrt{q}(z+z^{-1})}$. The integrals that enter
(\ref{Mnev}) and (\ref{Ahumada}) are related to each
other by the change of variable $s=\sqrt{q}(z+z^{-1})$.

\begin{acknowledgement}
We are indebted to Jessica Nelson and Wayne Lobb for some help
with the initial numerical experiments as well as for a number of
discussions, and to Markus Meringer for his help with
\texttt{GENREG}. The work of VE, JF and SL was supported, in part,
by the Australian Research Council Discovery grant DP0666632.
While this work was done, SL was at the University of South
Australia. The work of PZ was partially supported by the President
of Russian Federation grant NSh-U329.2006.1 and by the Russian
Foundation for Basic Research grant 05-01-00899.
\end{acknowledgement}

\end{document}